\documentclass[11pt,reqno]{amsproc}
\usepackage[semicolon,square,authoryear]{natbib}
 \usepackage{cite}
\usepackage[colorinlistoftodos, shadow]{todonotes}
\usepackage[T1]{fontenc}
\usepackage{amsmath,amscd,amssymb}
\usepackage{color}
\usepackage{graphicx}
\usepackage{psfrag,epsfig}
\usepackage{epstopdf}
\usepackage{multirow}
\usepackage{rotating}
\usepackage{fullpage}
\usepackage{subfigure}
\usepackage{sidecap}
\usepackage{enumerate}
\usepackage[small]{caption}
\usepackage{trimclip}
\usepackage{multirow}
\usepackage{listings}
\usepackage{longtable}
\usepackage[normalem]{ulem}

\numberwithin{equation}{section}

\usepackage[debug=false, colorlinks=true, pdfstartview=FitV, linkcolor=blue, citecolor=blue, urlcolor=blue]{hyperref}
\hypersetup{
    colorlinks=true,
    citecolor=blue,
    linkcolor=blue,
    filecolor=blue,      
    urlcolor=blue,
}
\usepackage{morefloats}
\newlength{\drop}
\definecolor{amethyst}{rgb}{0.6, 0.4, 0.8}
\definecolor{burgundy}{rgb}{0.5, 0.0, 0.13}

\title{On numerical stabilization in modeling double-diffusive viscous fingering} 
\author{\textbf{M.~Shabouei} and \textbf{K.~B.~Nakshatrala} \\ 
  {\small Department of Civil \& Environmental Engineering, 
    University of Houston \\
    Correspondence to: knakshatrala@uh.edu}}
\begin{document}
\def\AR{\clipbox{0pt 0pt .35em 0pt}{\textit{\bfseries A}}\kern-.05emR}

\begin{titlepage}
    \drop=0.1\textheight
    \centering
    \vspace*{\baselineskip}
    \rule{\textwidth}{1.6pt}\vspace*{-\baselineskip}\vspace*{2pt}
    \rule{\textwidth}{0.4pt}\\[\baselineskip]
         {\LARGE \textbf{\color{burgundy} 
		On numerical stabilization in modeling double-diffusive viscous fingering
         }}\\ [0.3\baselineskip]
    \rule{\textwidth}{0.4pt}\vspace*{-
   \baselineskip}\vspace{3.2pt}
    \rule{\textwidth}{1.6pt}\\[\baselineskip]
    \scshape
    An e-print of the paper is available on arXiv. \par
    
    \vspace*{0.5\baselineskip}
    Authored by \\[0.5\baselineskip]
    
    {\Large M.~Shabouei\par}
    {\itshape Graduate Student, University of Houston.}\\
     {\itshape Currently, Research Scientist, University of Notre Dame.}\\
    [\baselineskip]
    
    {\Large K.~B.~Nakshatrala\par}
    {\itshape Department of Civil \& Environmental Engineering \\
    University of Houston, Houston, Texas 77204--4003. \\ 
    \textbf{phone:} +1-713-743-4418, \textbf{e-mail:} knakshatrala@uh.edu \\
    \textbf{website:} http://www.cive.uh.edu/faculty/nakshatrala\par} 
    
\begin{figure*}[h]
\centering 
  \hspace{-0.30 cm}
  {\includegraphics[clip, scale=0.65]
  {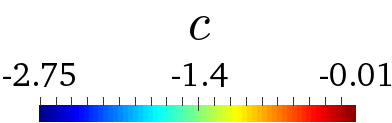}}
  	\vspace{0.05 cm}
  	\hspace{1.20 cm}
  {\includegraphics[clip, scale=0.65]
  {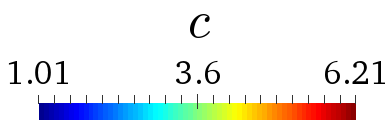}} \\
  	\vspace{0.05 cm}
  {\includegraphics[clip, scale=0.55]
  {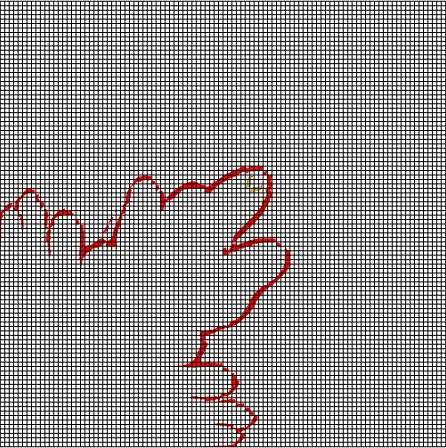}}
  	\hspace{0.35 cm}
  {\includegraphics[clip, scale=0.55]
  {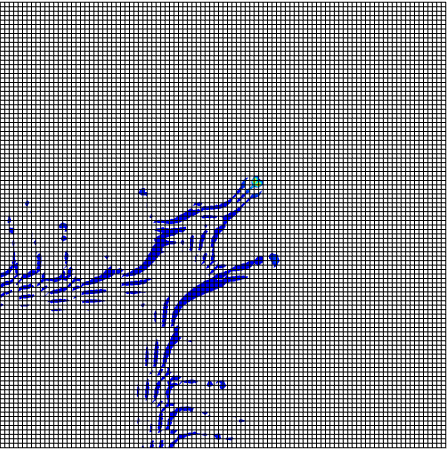}}
\end{figure*}
    {\small {\it This figure shows the unphysical concentration values got using the SUPG stabilized formulation on the quarter five-spot problem. The numerical solution for the concentration field violated the non-negative constraint (left) and the maximum principle (right). The mathematical model comprises coupled flow-thermal-transport equations.}}
    \vfill
    {\scshape 2019} \\
    {\small Computational \& Applied Mechanics 
    Laboratory} \par
\end{titlepage}

\begin{abstract}
A firm understanding and control of viscous fingering (VF) and miscible displacement will be vital to a wide range of industrial, environmental, and pharmaceutical applications, such as geological carbon-dioxide sequestration, enhanced oil recovery, and drug delivery.  We restrict our study to VF, a well-known hydrodynamic instability, in miscible fluid systems but consider double-diffusive (DD) effects---the combined effect of compositional changes because of solute transport and temperature. One often uses numerical formulations to study VF with DD effects. The primary aim of the current study is to show that popular formulations have limitations to study VF with DD effect. These limitations include exhibiting node-to-node spurious oscillations, violating physical constraints such as the non-negativity of the concentration field or mathematical principles such as the maximum principle, and suppressing physical instabilities. We will use several popular stabilized finite element formulations---the SUPG formulations and three modifications based on the SOLD approach---in our study. Using representative numerical results, we will illustrate two critical limitations. First, we document that these formulations do not respect the non-negative constraint and the maximum principle for the concentration field. We will also show the impact of these violations on how viscous fingers develop. Second, we show that these stabilized formulations, often used to suppress numerical instabilities, may also suppress physical instabilities, such as viscous fingering. Our study will be valuable to practitioners who use existing numerical formulations and to computational mathematicians who develop new formulations.
\end{abstract}
\maketitle
%

\section{INTRODUCTION AND MOTIVATION}
\label{Sec:S1_ViscFinger_Intro}
Hydrodynamic instabilities could occur in flows involving spatial disparity in the mobility ratio. Adverse mobility ratios arise because of variation in physical properties such as viscosity and density, and such disparities arise often in multi-phase fluid systems and non-isothermal systems. Two prevalent hydrodynamic instabilities are the Saffman-Taylor instability \citep{Saffman_1958_Hell_Shaw_VF}, also known as viscous fingering (VF), and the Rayleigh-Taylor instability \citep{Taylor_Density_instab_1950}, also known as density fingering. The latter instability arises because of a spatial disparity in the density, and the former, central to this paper, arises because of a spatial disparity in the viscosity. 

VF manifests itself in a wide range of industrial, environmental, chemical, and biological applications. These applications include enhanced oil recovery \citep{wooding_1976_VF_EOR},  geological carbon-dioxide sequestration \citep{Chen_1998_VF_2_JFM}, combustion \citep{zik_1998_VF_combust}, electrochemical processes \citep{deBruyn_1995_VF_electrochemical}, and chromatography \citep{Maes_2010_Exper_VF}. VF can occur in either immiscible or miscible systems. In an immiscible system of fluids, variation of surface tension at the interface of two fluids causes the fingering patterns. The capillary number governs instability in an immiscible system \citep{fernandez_2003_VF_immiscible}. Instability in a system with two miscible fluids arises because of a disorder in the flow field, and the velocity gradient dictates the fingering dynamics and mechanisms. In this paper, we shall restrict to miscible systems. 

VF and miscible displacement occur in porous media or Hele-Shaw cells when a less viscous fluid displaces a higher viscous fluid \citep{Saffman_1958_Hell_Shaw_VF}. Because of the faster movement of a low viscous fluid, any small disturbance or perturbation to the interface of miscible fluids generates finger-like structures. Heterogeneity in material properties, such as permeability and diffusivity, or nonlinearities can induce these perturbations. Two primary sources of nonlinearity are fluid properties depending on the concentration or temperature, and nonlinear chemical reactions. 

Viscosity contrasts, which drive VF, can occur because of compositional changes \citep{tan_homsy_1988_VF}, temperature contrasts \citep{mishra_2010_DD_VF}, or chemical reactions \citep{hejazi_2010_VF}. The initial works have studied the so-called classical VF, which arises because of compositional changes. VF occurring because of the combined effect of composition changes and temperature contrasts falls under the study of double-diffusive (DD) VF. However, DD effects can occur for factors other than temperature. Different solutes or mass components with different diffusive rates can also trigger VF, and such scenarios also form a part of studies on DD effects \citep{mishra_2010_DD_VF}. In this paper, we restrict VF with DD effects because of temperature and single component solute transport. 

One can study VF with DD effects by considering miscible displacement between fluids in non-isothermal conditions \citep{chan_1997_DF_thermal}. However, it is essential to note that the dynamics of VF with DD effects can differ from that of the classical VF. For example, the stable case under the classical VF---a lower viscosity fluid displacing another fluid with a higher viscosity---may not be stable because of temperature contrasts. Since the time scale of heat transfer is faster than that of mass transfer, DD effects can influence the stability characteristics of VF. There have been several theoretical investigations of the classical VF and DD effects. Some notable ones on the classical VF include \citep{Saffman_1958_Hell_Shaw_VF, Raghavan_1971_VF,Saffman_1986_JFM_VF, Chang_1986_VF, tan_homsy_1986_VF, riaz_2004_VF} and \citep{hejazi_2010_VF}. Some representative works on VF with DD effects are \citep{chan_1997_DF_thermal,pritchard_2004_DD, pritchard_2009_DD_VF, nagatsu_2009_experimental_DD, 
islam_2010_DD_VF_part_1,islam_2010_DD_VF_part_2}. However, it is not possible to get analytical solutions for problems involving VF, as the mathematical models involve a system of nonlinear, coupled partial differential equations. Therefore, many studies have used numerical simulations to gain insight into the impact of material and physical properties, including viscosity, density, permeability, and diffusivity dispersion, on hydrodynamic instabilities. For example, see \citep{tan_homsy_1988_VF, tan_homsy_1992_VF, zimmerman_homsy_1992_VF, zimmerman_homsy_1992_VF_2, tchelepi_1993_VF, Chen_1998_VF_2_JFM, riaz_2006_VF}. Researchers have also been developing several numerical methods to facilitate qualitative and quantities studies of VF. For example, see \citep{scovazzi_2013_JCP_VF,chang2017variational,Li_Riviere_2015numerical}. However, these works have addressed the classical VF (i.e., with single diffusivity).

Given its importance of numerical simulations in the study on VF, especially with double-diffusive effects, we will undertake a critical study on how some popular numerical formulations perform and document their limitations. We will consider a mathematical model capable of simulating VF with double-diffusive effects by coupling flow, transport, and thermal processes. We shall show that some popular numerical stabilized formulations (e.g., SUPG and SOLD), used to avoid spurious numerical artifacts (e.g., node-to-node spurious oscillations) for advection-dominated problems, may suppress mathematical and physical VF instabilities. Many of the existing numerical formulations can give unphysical negative values for the concentration field. These unphysical values will adversely affect the accuracy of the numerical simulators in predicting physical instabilities like VF. We document these formulations violate the non-negative constraint in the concentration field and show the effects on VF profiles. The benefit of our study to the scientific community will be twofold. First, the results we present in this paper will provide valuable information to practitioners and forewarn them on the potential problems in using numerical solutions of problems involving VF. Second, our study will guide computational mathematicians in developing new accurate numerical formulations for VF with double-diffusive effects. 

We organize the rest of this paper as follows. Section \ref{Sec:S2_ViscFinger_GE} presents a mathematical model for double-diffusive effects by coupling flow, transport, and thermal processes. In Section \ref{Sec:S3_Stab_FEM}, we outline the stabilized formulations used in our study. In Section \ref{Sec:S4_Stab_NR}, we illustrate using numerical results on how several popular numerical stabilizations perform and document their limitations. We will end the paper with concluding remarks in Section \ref{Sec:S5_Stab_CR}.

\section{GOVERNING EQUATIONS:~COUPLED FLOW, TRANSPORT, AND THERMAL PROCESSES}
\label{Sec:S2_ViscFinger_GE}
We consider the flow of an incompressible single-phase miscible fluid in rigid porous media, including transport and thermal effects. We now present a mathematical model in the form of a system of coupled partial differential equations (PDEs), each of which is elliptic or parabolic. The governing equations for the transport and thermal subproblems are transient, while the governing equations for the flow subproblem are quasistatic. The mathematical model involves a two-way coupling: the viscosity of the fluid is a (nonlinear) function of the concentration and temperature, and the transport and thermal equations depend on the flow velocity of the fluid. 

Let $\Omega \subset \mathbb{R}^{d}$ be an open bounded domain, where ``$d$'' denotes the number of spatial dimensions. Let $\partial \Omega := \overline{\Omega} - \Omega$ denote the boundary of the domain, where $\overline{\Omega}$ denotes
the set closure of $\Omega$. A spatial point is denoted by $\mathbf{x} \in \overline{\Omega}$. The time is denoted by $t \in ]0,\mathcal{I}[$, where $\mathcal{I}$ denotes the total time of interest. The gradient and divergence operators with respect to $\mathbf{x}$ are, respectively, denoted by $\mathrm{grad}[\cdot]$ and $\mathrm{div} [\cdot]$. The governing equations for the \emph{flow} subproblem, governed by Darcy equations, are written as follows:
\begin{subequations}
  \begin{align}
    \label{Eqn:Darcy_GE}
    & \frac{\mu(c,\theta)}{k}  
    \mathbf{v}(\mathbf{x},t) + \mathrm{grad}
    [p(\mathbf{x},t)] = \mathbf{0} \quad 
    && \mathrm{in} \; \Omega \times ]0,\mathcal{I}[ \\
    \label{Eqn:Darcy_IncompCon}
    & 
    \mathrm{div}[\mathbf{v}(\mathbf{x},t)] 
    = \varphi(\mathbf{x},t) \quad 
    && \mathrm{in} \; \Omega \times ]0,\mathcal{I}[ \\
    \label{Eqn:Darcy_PressureBC}
    & p(\mathbf{x},t) = p_0(\mathbf{x},t) \quad 
    && \mathrm{on} \; \Gamma^{p} \times ]0,\mathcal{I}[ \\
    \label{Eqn:Darcy_VelocityBC}
    &\mathbf{v}(\mathbf{x},t) \bullet \widehat{\mathbf{n}}
    (\mathbf{x}) = v_n(\mathbf{x},t) \quad 
    && \mathrm{on} \; \Gamma^{v} \times ]0, \mathcal{I}[
  \end{align}
\end{subequations}
where $\mu(c,\theta)$ is the viscosity of the fluid, $\mathbf{v}(\mathbf{x},t)$ is the velocity vector field, and the pressure field is denoted by $p(\mathbf{x},t)$. $k$ is the permeability and $\varphi(\mathbf{x},t)$ denotes the mass production. $p_0(\mathbf{x},t)$ is the prescribed pressure on the boundary and $v_n(\mathbf{x},t)$ is the prescribed normal component of the velocity vector field on the boundary. $\Gamma^{p}$ and $\Gamma^{v}$ are, respectively, the boundaries of the domain on which the pressure and normal component of the velocity are prescribed. It is assumed that the velocity and its gradient are small so that the inertial effects, and consequently the convective term, can be neglected in the balance of linear momentum.

The \emph{transport} of chemical species is modeled using advection-diffusion-reaction (ADR) equations, and the corresponding governing equations take the following form:
\begin{subequations}
  \begin{align}
    \label{Eqn:ADRs_for_i}
    &\frac{\partial c(\mathbf{x},t)}{\partial t} 
    + \mathrm{div} \left[ \mathbf{v} c - d_m \; 
    \mathrm{grad}[c] \right] = 
    r(\mathbf{x},t,c) \quad 
    &&\mathrm{in} \; \Omega \times ]0, 
    \mathcal{I}[ \\
    \label{Eqn:ADRs_Dirchlet_i}
    &c(\mathbf{x},t) = c^{\mathrm{p}}(\mathbf{x},t) 
    \quad 
    &&\mathrm{on} \; \Gamma^{c} \times ]0, 
    \mathcal{I}[ \\
    \label{Eqn:ADRs_Neumann_i}
    &\left( \mathbf{v} c - d_m \; \mathrm{grad} [c]
    \right) \bullet \widehat{\mathbf{n}}
    (\mathbf{x}) = h^{\mathrm{p}}(\mathbf{x},t) \quad 
    &&\mathrm{on} \; \Gamma^{h} \times ]0, \mathcal{I}[ \\
    \label{Eqn:ADRs_IC_i}
    &c(\mathbf{x},t=0) = c^{0}(\mathbf{x}) \quad 
    &&\mathrm{in} \; \Omega
  \end{align}
\end{subequations}
where $c(\mathbf{x},t)$ is the concentration of the chemical species and $d_m$ is the corresponding diffusivity. $r(\mathbf{x},t,c)$ is the reactive component of the volumetric source. $c^{0}(\mathbf{x})$ is the initial condition, $c^{\mathrm{p}}(\mathbf{x},t)$ is the prescribed concentration on the boundary and $h^{\mathrm{p}}(\mathbf{x},t)$ is the prescribed flux on the boundary. $\Gamma^{c}$ and $\Gamma^{h}$, respectively, denote the Dirichlet and Neumann parts of the boundary for the transport subproblem. We assume that the ADR equations are in mass fraction form, and hence we have $0 \leq c \leq 1$. 

The governing equations for the \emph{thermal} subproblem take the following form:
\begin{subequations}
  \begin{align}
    \label{Eqn:Temp_equ}
    &\frac{\partial \theta(\mathbf{x},t)}{\partial t} 
    + \mathrm{div}[\mathbf{v} 
    \theta - \kappa_{\theta} \, 
    \mathrm{grad}[\theta]] = 
    q(\mathbf{x},t,\theta) \quad 
    &&\mathrm{in} \; \Omega \times ]0, \mathcal{I}[ \\
    \label{Eqn:Temp_Dirchlet}
    &\theta(\mathbf{x},t) = \theta^{\mathrm{p}}(\mathbf{x},t) 
    \quad 
    &&\mathrm{on} \; \Gamma^{\theta} \times ]0, 
    \mathcal{I}[ \\
    \label{Eqn:Temp_Neumann}
    &\left(\mathbf{v} \theta - \kappa_{\theta} \mathrm{grad}
    [\theta] \right) \bullet \widehat{\mathbf{n}}(\mathbf{x}) = 
    s^{\mathrm{p}}(\mathbf{x},t) \quad 
    &&\mathrm{on} \; \Gamma^{s} 
    \times ]0, \mathcal{I}[ \\
    \label{Eqn:Temp_IC}
    &\theta(\mathbf{x},t=0) = \theta^{0}(\mathbf{x}) \quad 
    &&\mathrm{in} \; \Omega
  \end{align}
\end{subequations}
where $\theta(\mathbf{x},t)$ is the temperature and $\kappa_{\theta}$ is the thermal diffusivity. $q(\mathbf{x},t,\theta)$ is the thermal volumetric source. $\theta^{0}(\mathbf{x})$ is the prescribed initial temperature. $\theta^{\mathrm{p}}(\mathbf{x},t)$ and $s^{\mathrm{p}}(\mathbf{x},t)$ are, respectively, the prescribed temperature and the prescribed heat flux on the boundary. $\Gamma^{\theta}$ and $\Gamma^{s}$, respectively, denote the Dirichlet and Neumann parts of the boundary for the thermal subproblem.

For simplicity, we assume that the medium properties are isotropic, implying that permeability $k$, mass diffusivity $d_m$, and thermal diffusivity $\kappa_{\theta}$ are all scalars. 

\section{STABILIZED FINITE ELEMENT FORMULATIONS}
\label{Sec:S3_Stab_FEM}
Given the mathematical complexity of the governing equations, it is not possible always to find analytical solutions to these equations. One often seeks numerical solutions to understand the underlying dynamics. For concreteness, we resort to the finite element method in this paper. However, the conclusions will be valid even for other discretization methods like the finite difference method, finite volume method, and lattice Boltzmann method. 

Two different numerical instabilities could arise in the numerical solutions of advection-diffusion (AD) and advection-diffusion-reaction (ADR) systems. The first instability is spurious node-to-node oscillations, which often manifest in a solution field for advection-dominated problems (i.e., problems with high P\'eclet numbers). Problems with strong advection exhibit steep gradients in the solution field, which manifest as interior and boundary layers. A numerical formulation may not resolve such sharp features on coarse computational grids. For example, it is a well-known fact the standard single-field Galerkin finite element formulation provides inaccurate numerical solutions (polluted with node-to-node oscillations) when the element P\'eclet number exceeds unity \citep{Gresho_Sani_v1}. (The element P\'eclet number is a convenient non-dimensional number in terms of advection velocity, diffusivity, and the characteristic element size.) In principle, it is always possible to bring down the element P\'eclet number by making the mesh finer. However, such a procedure is not attractive because of the computational cost. One can instead use stabilized formulations, which suppress spurious numerical oscillations \citep{Augustin_2011_Stabilizers}. Such stabilized formulations will form the focus of this paper. 

The second instability pertains to numerical solutions violating physical constraints, such as the non-negative constraint, or mathematical principles, such as the maximum principle. The governing equation for the transport problem (i.e., equation \eqref{Eqn:ADRs_for_i}) is a parabolic PDE and such a PDE possesses several important mathematical properties, including the maximum principle \citep{Gilbarg_Trudinger}. Physics dictates that quantities like concentration and (absolute) temperature can attain only non-negative values. It is now well-documented that many popular numerical methods, such as the finite element method \citep{Liska_2008_FEM_non_negative,Nagarajan_Nakshatrala_IJNMF_2010, 
Nakshatrala_2013_ADR_non_negative, Nakshatrala_2013_non_negative}, finite volume method \citep{Le_Potier_2005_FVM_non_negative}, and lattice Boltzmann method \citep{karimi_2015_LBM_non_negative}, violate the maximum principle and the non-negative constraint for diffusion-type equations, especially when the diffusivity is strongly anisotropic. \emph{In this paper, we will show that numerical formulations violate the maximum principle and the non-negative constraint even under isotropic diffusivity for VF with DD effects.}

It is important not to confuse numerical instabilities with that of hydrodynamic instabilities, which are physical, such as viscous fingering. It is desirable that a numerical formulation is devoid of numerical instabilities but able to capture physical instabilities. We discuss below and illustrate using numerical results in the next section that some popular formulations do not achieve this desired feature.

\subsection{Stabilized formulations for ADR equations}
Herein, we shall utilize the Streamline Upwind Petrov-Galerkin (SUPG) formulation \citep{Brooks_Hughes_CMAME_1982_v32_p199} and three modifications of SUPG based on the Spurious Oscillations at Layers Diminishing (SOLD) approach \citep{Hughes_SUPG_modified_1986,Johnson_Crosswind_1987,2007_John_Knobloch_CMAME_v196_p2197_p2215}. These formulations are among the most popular stabilized FEM formulations for solving AD and ADR systems, similar to the one given by equations \eqref{Eqn:ADRs_for_i}--\eqref{Eqn:ADRs_IC_i}. 

We introduce the following function spaces:
\begin{subequations}
  \begin{align}
    \label{Eqn:Function_Space_C}
    \mathcal{C}_t &:= \left\{c(\mathbf{x}, 
    t) \in H^{1}(\Omega) \; \big| 
    \; c(\mathbf{x},t) = c^{\mathrm{p}}
    (\mathbf{x},t) \; \mathrm{on} \; 
    \Gamma^{c} \times ]0, \mathcal{I}[ 
    \right\} \\
    \label{Eqn:Function_Space_W}
    \mathcal{W} &:= \left\{w(\mathbf{x}) 
    \in H^{1}(\Omega) \; \big| 
    \; w(\mathbf{x}) = 0 \; \mathrm{on} 
    \; \Gamma^{c} \right\} 
  \end{align}
\end{subequations}
where $H^{1}(\Omega)$ is a standard Sobolev space \citep{Evans_PDE}. The standard $L_2$ inner-product for given two fields $a(\mathbf{x},t)$ and $b(\mathbf{x},t)$ over a set $\mathcal{B}$ is defined as follows:
\begin{align}
  \label{Eqn:L2_Inner_Product}
  \left(a;b\right)_{\mathcal{B}} = \int 
  \limits_{\mathcal{B}} a(\mathbf{x}) 
  \bullet b(\mathbf{x}) \; \mathrm{d} 
  \mathcal{B}
\end{align}
The subscript on the inner-product will be dropped if $\mathcal{B} = \Omega$. 

\subsubsection{The SUPG formulation} 
Find $c(\mathbf{x},t) \in \mathcal{C}_t$ such that we have
\begin{align}
  \label{Eqn:SUPG_Formulation}
  \left(w; \frac{\partial c}{\partial t}\right) 
  &+ (\mathrm{grad}[w] \bullet \mathbf{v};c) 
  + (\mathrm{grad}[w];\mathbf{D}(\mathbf{x},t) 
  \mathrm{grad}[c]) \nonumber \\ 
  &+ \sum_{e = 1}^{Nele} \left(\tau \, 
  \mathbf{v} \bullet \mathrm{grad}[w]; 
  \frac{\partial c}{\partial t} + \mathrm{div} 
  \left[\mathbf{v}c - \mathbf{D}(\mathbf{x},t) 
  \mathrm{grad}[c] \right] - f - r 
  \right)_{\Omega_e} \nonumber \\
  &= (w;f) + \left(w;h^{\mathrm{p}} 
  \right)_{\Gamma^{h}} \quad \forall 
  w(\mathbf{x}) \in \mathcal{W}
\end{align}
where $Nele$ denotes the number elements in the finite element mesh, $\tau$ is the stabilization parameter, and $\Omega^{e}$ denotes an element in the finite element mesh. Note that $\bar{\Omega} = \bigcup_{e = 1}^{Nele} \bar{\Omega}^{e}$, where a superposed bar denotes the set closure. We use the stabilization parameter proposed by \citet{2007_John_Knobloch_CMAME_v196_p2197_p2215}, with the following form: 
\begin{align}
  \label{Eqn:StabParam_SUPG}
  \tau (\mathbf{v}) = \frac{h_{\Omega_e}}{2 
  \| \mathbf{v} \|} \xi_0 \left(\mathbb{P}\mathrm{e}_{h}\right)
\end{align}
where $h_{\Omega_e}$ is the maximum element length, $\xi_0$ is the upwind function, and $\mathbb{P} \mathrm{e}_{h}$ is the local (element) P\'eclet number. The definitions for $\xi_0$ and $\mathbb{P}\mathrm{e}_h$ are as follows: 
\begin{align}
  &\xi_{0} \left(\chi \right) = \coth \left( \chi\right) - \frac{1}{\chi} \\
  &\mathbb{P}\mathrm{e}_{h} = \frac{h_{\Omega_e} 
    \| \mathbf{v} \|}{2 \lambda_{\mathrm{min}}}
\end{align}
where $\lambda_{\mathrm{min}}$ denotes the minimum eigenvalue of the diffusivity tensor, and $\|\cdot\|$ denotes the Euclidean norm. 

\subsubsection{SOLD modification with isotropic artificial diffusion} 
To diminish the oscillations in the numerical solutions under the SUPG formulation, \citet{Hughes_SUPG_modified_1986} proposed a modification by adding artificial diffusion along the upwind direction. The resulting weak form can be obtained by augmenting the following term to the left hand side of equation \eqref{Eqn:SUPG_Formulation}:
\begin{align}
	\label{Eqn:SOLD_iso_Formulation}
	\sum_{e = 1}^{Nele} \left(\tau_1 \, 
	\mathbf{v}^{\|} 
	\bullet \mathrm{grad}[w]; 
	\frac{\partial c}{\partial t} 
	+ \mathrm{div} \left[\mathbf{v}c 
	- \mathbf{D}(\mathbf{x},t) 
	\mathrm{grad}[c] \right] - f 
	- r \right)_{\Omega_e}
\end{align}
where $\mathbf{v}^{\|}$ is defined as follows \citep{2007_John_Knobloch_CMAME_v196_p2197_p2215}:
\begin{align*}
	\mathbf{v}^{\|} = \left\{ \begin{array}{lcl}
	\frac{(\mathbf{v} \bullet \mathrm{grad}[c]) 
	\mathrm{grad}[c]}{ \|\mathrm{grad}[c]\|^2 }
	& \mbox{for} & \mathrm{grad}[c] \neq 
	\mathbf{0} \\ 
	\mathbf{0} & \mbox{for} & \mathrm{grad}[c] 
	= \mathbf{0} 
	\end{array}\right.
\end{align*}
and the stabilization parameter is given by the following expression \citep{Hughes_SUPG_modified_1986}:
\begin{align*}
  \tau_1 = \mathrm{max} \{ 0, \tau (\mathbf{v}^{\|}) - \tau (\mathbf{v}) \}
\end{align*}
%
\subsubsection{SOLD modification with cross-wind artificial diffusion} 
\citet{Johnson_Crosswind_1987} proposed an alternate modification by adding artificial diffusion in the cross-wind direction to the SUPG formulation. They defined the following projection operator, which maps a vector onto a plane orthogonal to $\mathbf{v}$:
\begin{align*}
	\mathbb{P}^{\perp} = \left\{ \begin{array}{lcl}
	\mathbf{I} - \frac{\mathbf{v} \otimes 
	\mathbf{v}}	{ \|\mathbf{v}\|^2 }
	& \mbox{for} & \mathbf{v} \neq \mathbf{0} \\
	\mathbf{0} & \mbox{for} & \mathbf{v} 
	= \mathbf{0} 
	\end{array}\right.
\end{align*}
where $\mathbf{I}$ is the identity tensor. The resulting weak form under this modification is obtained by adding the following term to the left hand side of equation \eqref{Eqn:SUPG_Formulation}:
\begin{align}
	\label{Eqn:SOLD_cross_Formulation}
	\sum_{e = 1}^{Nele} \left( \tau_2 \, 
	\mathbb{P}^{\perp} \; \mathrm{grad}[w]; 
	\mathrm{grad}[c] \right)_{\Omega_e}
\end{align}
where the stabilization parameter $\tau_2$ is defined as follows: 
\begin{align*}
	\tau_2 = \mathrm{max} 
	\{ 0, \| \mathbf{v} \| \; 
	h_{\Omega_e}^{2/3} 
	- \lambda_{\mathrm{min}} \}
\end{align*}

We also use a third modification to the SUPG formulation, obtained by incorporating both isotropic diffusion and cross-wind diffusion. The resulting weak form is obtained by adding the artificial diffusion terms \eqref{Eqn:SOLD_iso_Formulation} and \eqref{Eqn:SOLD_cross_Formulation} to the left hand side of equation \eqref{Eqn:SUPG_Formulation}.

We will use either the SUPG formulation or one of the modifications for solving the transport subproblem, and use the SUPG formulation for the thermal subproblem.

\subsection{Mixed weak formulation for Darcy equations}
For completeness, we provide the details of the weak formulation for solving the flow equations, even though this formulation is not a stabilized formulation. We used the classical mixed formulation, also known as the Galerkin mixed formulation, to solve the Darcy equations. However, it is important to note that all combinations of interpolations for the fields variables---the velocity and the pressure---give stable results. The celebrated Ladyzhenskaya-Babu\v ska-Brezzi (LBB) \emph{inf-sup} stability condition provides the mathematical theory to address the stability of a mixed formulation \citep{Brezzi_Fortin}. In our study, we employ quadrilateral elements with second-order interpolation (Q9) for the velocity field and first-order interpolation (Q4) for the pressure field. This combination of the interpolations for quadrilateral elements satisfies the LBB stability condition, as any $Q_{n}Q_{n-1}$ combination of interpolations satisfies the LBB stability condition \citep{Brezzi_Fortin}. Our choice ensures that there will be no spurious oscillations in the pressure field. 

For a fixed time-level $t$, the classical mixed formulation for equations 
\eqref{Eqn:Darcy_GE}--\eqref{Eqn:Darcy_VelocityBC} reads:~Find 
$\mathbf{v}(\mathbf{x},t)\in\mathcal{V}$ and $p(\mathbf{x},t) \in L_2(\Omega)$ such that we have
\begin{align}
  \label{Eqn:Dr_mixed_Formulation}
  (\mathbf{w}; \mu \mathbf{K}^{-1} \mathbf{v}) 
  - (\mathrm{div}[\mathbf{w}]; p) 
  - (q; \mathrm{div}[\mathbf{v}])
  = (\mathbf{w}; \rho \mathbf{b}) - (\mathbf{w} \bullet \widehat{\mathbf{n}}
  ; p_0)_{\Gamma^p} 
  \quad \forall \mathbf{w}(\mathbf{x}) \in \widetilde{\mathcal{W}}, q(\mathbf{x}) \in L_2(\Omega)
\end{align}
wherein the following function spaces are used: 
\begin{subequations}
  \begin{align}
    \label{Eqn:Function_Space_V}
    \mathcal{V} &:= \left\{\mathbf{v}(\mathbf{x}, 
    t) \in (L_{2}(\Omega))^{nd} 
    \; \big| \; \mathrm{div}[\mathbf{v}] 
    \in L_{2}(\Omega), \; \mathbf{v}(\mathbf{x},t) 
    \bullet \widehat{\mathbf{n}}
    (\mathbf{x}) = v_n(\mathbf{x},t) \; \mathrm{on} 
    \; \Gamma^{v} \times ]0, \mathcal{I}[ 
    \right\} \\
    \label{Eqn:Function_Space_U}
    \widetilde{\mathcal{W}} &:= \left\{\mathbf{w}(\mathbf{x}) 
    \in (L_{2}(\Omega))^{nd} \; \big| 
    \; \mathrm{div}[\mathbf{w}] 
    \in L_{2}(\Omega), \; \mathbf{w}(\mathbf{x}) 
    \bullet \widehat{\mathbf{n}}
    (\mathbf{x}) = 0 \; \mathrm{on} \; \Gamma^{v} 
    \right\}
  \end{align}
\end{subequations}

\section{ON NUMERICAL SOLUTIONS USING STABILIZED FORMULATIONS}
\label{Sec:S4_Stab_NR}
In this section, we highlight two potential problems in numerical modeling of VF and, in particular, when we deal with double-diffusive effects. First, numerical formulations for transport equations may render unphysical solutions, polluted with negative values in the concentration field. These numerical solutions may violate maximum principles, making these solutions unreliable and unphysical. Second, if one uses a stabilized formulation to suppress numerical instabilities, if not judicious in the selection, the formulation may also suppress underlying physical instabilities like the viscous fingers.

\begin{figure}
\centering 
  \includegraphics[clip, scale=0.40]
  {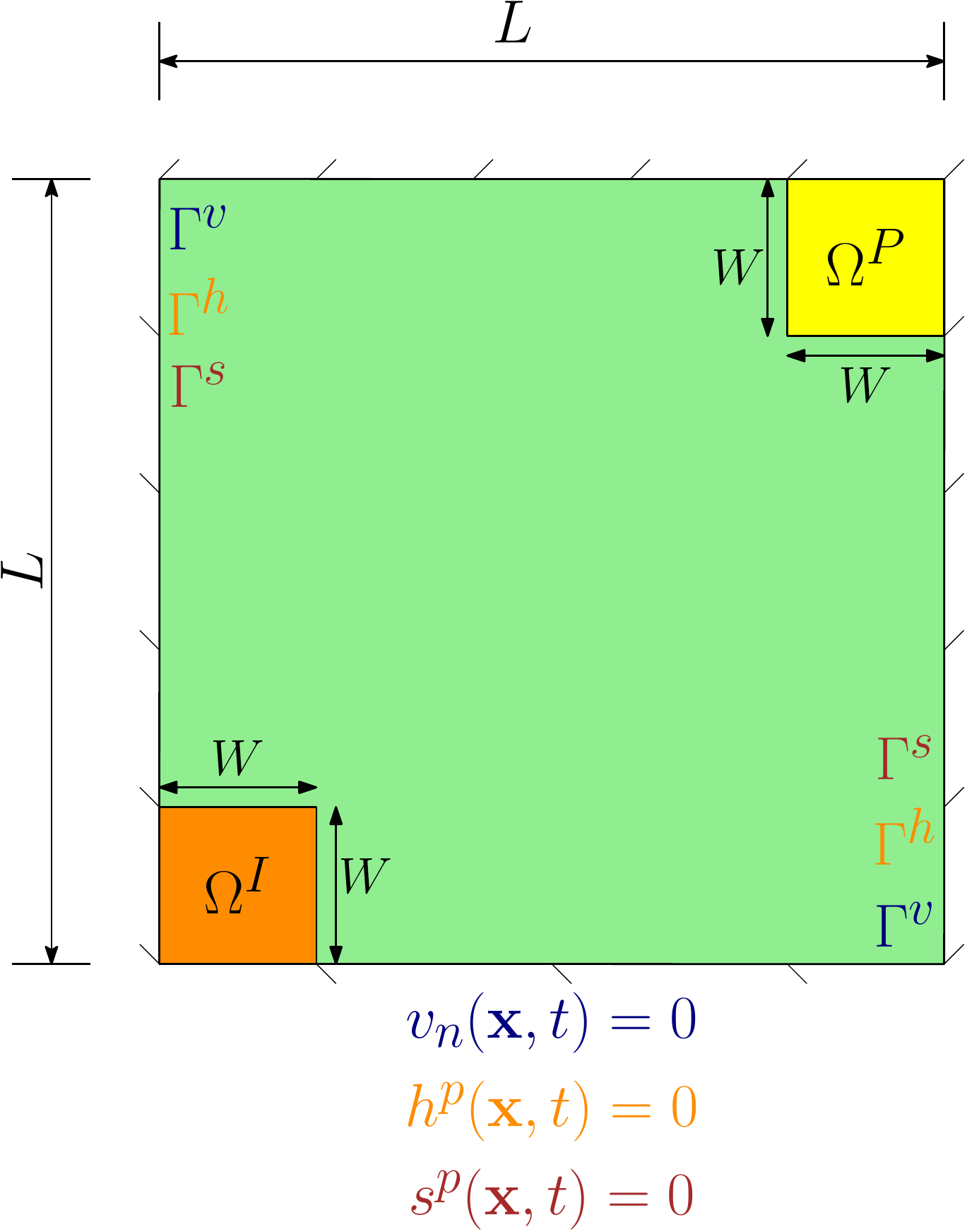}
  \caption{\textsf{Quarter five-spot problem:}~Pictorial description of the problem including initial and boundary conditions. The computational domain includes a unit square ($L = 1$) and two small squares with length $W$. The bottom left corner square ($\Omega^I$) is injection well and the top right one ($\Omega^P$) is production well. \label{Fig:pict_descript_case_II}}
\end{figure}

We alter the quarter five-spot problem, used often in studies on flow through porous media, to fit the studies on double-diffusive effects. Figure \ref{Fig:pict_descript_case_II} provides a pictorial description of the problem. The computational domain is a square with length $L = 1$. There are also two small squares of length $W = L/10$ near the bottom left, the injection well $\Omega^I$, and top right, the production well $\Omega^P$, corners.
One needs to introduce small perturbations into the prescribed data (e.g., initial condition, medium properties like permeability, or boundary conditions) to trigger physical instabilities like the viscous fingers. Instead, one can introduce a chemical reaction into transport equations, resulting in advection-diffusion-reaction (ADR) equations, to trigger physical instabilities. In the studies presented in this section, we introduce a chemical reaction in both transport and thermal subproblems near the production well ($\Omega^{P}$) to trigger viscous fingers.

From the injection well of $\Omega^I$, displacing fluid mass ($\varphi^I$) is injected into the domain which is already filled by displaced fluid at rest. The top right square ($\Omega^P$) is production well that has mass sink ($\varphi^P$). Elsewhere, $\varphi^I = \varphi^P = 0$. Homogeneous velocity is enforced 
on the entire boundary (i.e., $v_n(\mathbf{x},t) = 0 \; \mathrm{on} \; \partial \Omega \times ]0, \mathcal{I}[$). The zero 
fluxes are prescribed on the entire boundary 
for ADR equations (i.e., $h^{\mathrm{p}}(\mathbf{x},t) = 0 
\; \mathrm{on} \; \partial \Omega \times ]0, \mathcal{I}[$) and heat equations (i.e., 
$s^{\mathrm{p}}(\mathbf{x},t) = 0 \; \mathrm{on} \; \partial \Omega \times ]0, \mathcal{I}[$).
As before, the viscosity depends exponentially on the concentration and temperature. That is, $\mu = \mu_0 \exp\left[ R_c (1-c) + R_{\theta} (1-\theta) \right]$. The injection source and production sink are $\varphi^I = \varphi^P = 0.1$, respectively. The maximum principle implies that the concentration field should lie between 0 and 1; that is, $0 \leq c(\mathbf{x},t) \leq 1$.

\begin{table}
  \begin{center}
    \caption{\textsf{Quarter five-spot:}~Parameters 
    used in the problem.
      \label{Tab:Input_case_II_dispersion}}
    \begin{tabular}{c c}
      \hline
      Parameter                 & Value \\ \hline\hline
      $d_m$			& $10^{-7}$ \\
      $f(\mathbf{x},t)$		& $\varphi^I \;\; \mathrm{in} 
      \;\; \Omega^I$ \\
      $r(\mathbf{x},t,c)$	& $-\varphi^Pc \;\; \mathrm{in} 
      \;\; \Omega^P$ \\
      $g(\mathbf{x},t)$		& $\varphi^I \;\; \mathrm{in} 
      \;\; \Omega^I$ \\
      $q(\mathbf{x},t,\theta)$	& $-\varphi^P \theta 
      \;\; \mathrm{in} \;\; \Omega^P$ \\
      $R_c$ & 2 \\
      $R_{\theta}$ & 2 \\ 
      $\mu_0$ & 1 \\
      \hline 
    \end{tabular}
  \end{center} 
\end{table}

\begin{figure}
\centering  
  {\includegraphics[clip, scale=0.5]
  {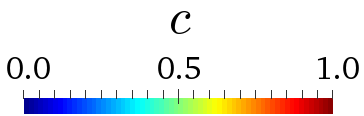}} \\
  %
  %
  \subfigure[SUPG at $t = 100$]
  {\includegraphics[clip, scale=0.3]  
  {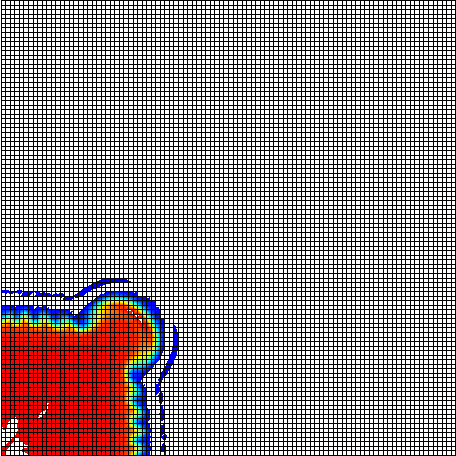}}
  \hspace{0.15 cm}
  \subfigure[SUPG at $t = 175$]
  {\includegraphics[clip, scale=0.3]
  {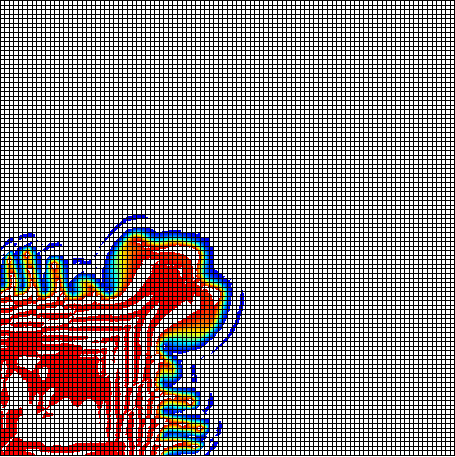}}
  \hspace{0.15 cm}
  \subfigure[SUPG at $t = 250$]
  {\includegraphics[clip, scale=0.3]
  {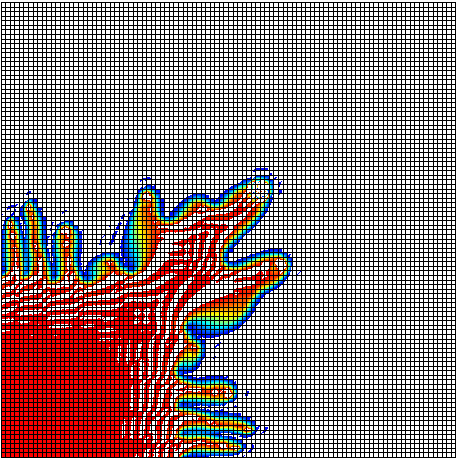}} \\  
  \vspace{0.015 cm}
  %
  \subfigure[SUPG + isotropic diffusion at $t = 100$]
  {\includegraphics[clip, scale=0.3]  
  {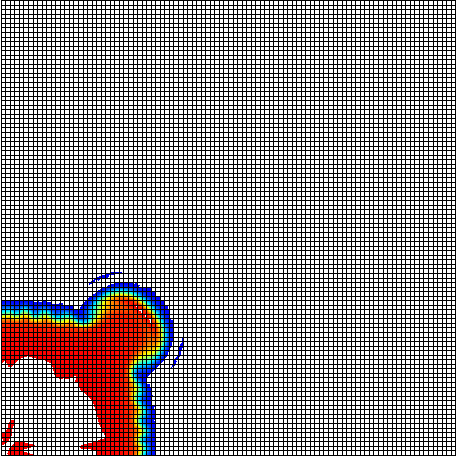}}
  \hspace{0.15 cm}
  \subfigure[SUPG + isotropic diffusion at $t = 175$]
  {\includegraphics[clip, scale=0.3]
  {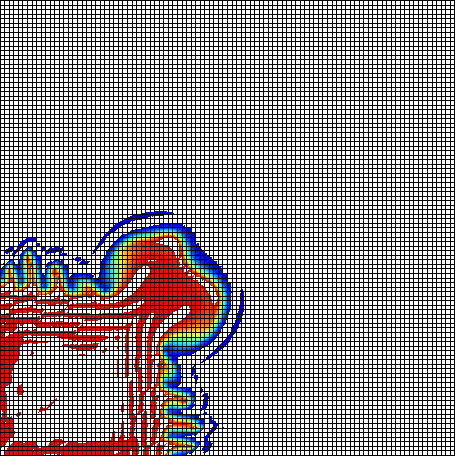}}
  \hspace{0.15 cm}
  \subfigure[SUPG + isotropic diffusion at $t = 250$]
  {\includegraphics[clip, scale=0.3]
  {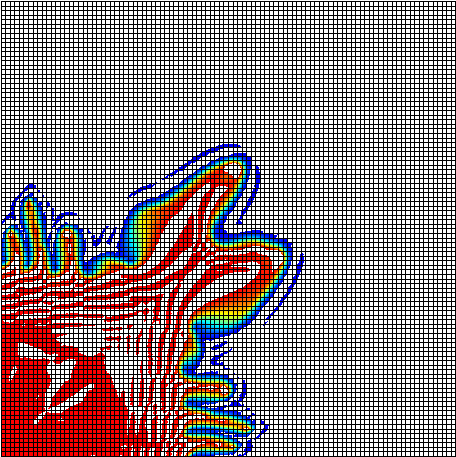}} \\  
  \vspace{-0.04 cm}
  %
  \subfigure[SUPG + crosswind diffusion at $t = 100$.]
  {\includegraphics[clip, scale=0.3]  
  {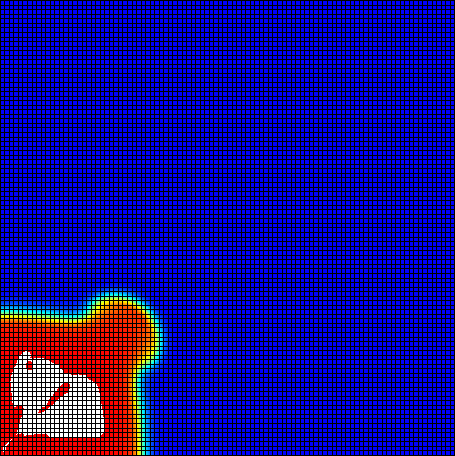}}
  \hspace{0.15 cm}
  \subfigure[SUPG + crosswind diffusion at $t = 175$]
  {\includegraphics[clip, scale=0.3]
  {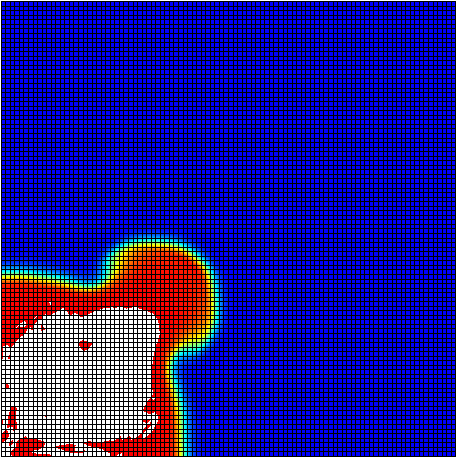}}
  \hspace{0.15 cm}
  \subfigure[SUPG + crosswind diffusion at $t = 250$]
  {\includegraphics[clip, scale=0.3]
  {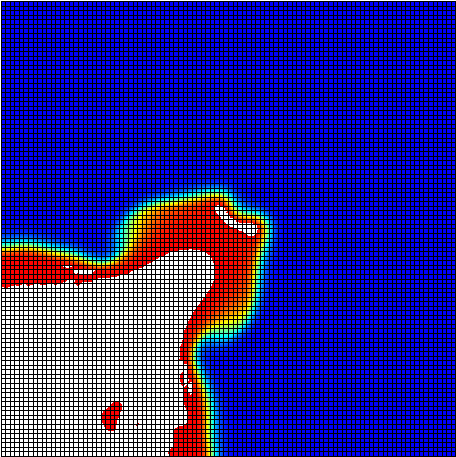}}  
  \caption{\textsf{Violation of the maximum principle:}~These figures show concentration  
  at different time steps for coupled flow, 
  ADR and temperature equations under numerical 
  stabilizers. The upper figures are for SUPG 
  formulation (equation \eqref{Eqn:SUPG_Formulation}). 
  the middle figures are for modified SUPG in which 
  we added isotropic artificial diffusion (equation 
  \eqref{Eqn:SOLD_iso_Formulation}) and bottom figures 
  are for modified SUPG with crosswind artificial 
  diffusion (equation \eqref{Eqn:SOLD_cross_Formulation}). 
  The results are obtained for parameters provided 
  in Table \ref{Tab:Input_case_II_dispersion}. 
  The white regions represent the regions 
  in which the concentration has violated 
  the maximum principle.
  The injection source and production 
  sink are $\varphi^I = \varphi^P = 0.1$.
  The mesh used for this problem is $100 \times 
  100$ structured quadrilateral elements.
  The numerical results clearly violate the 
  \emph{maximum principle} (and in consequence 
  non-negative constraint) for the concentration. 
  \label{Fig:conc_case_I}}
\end{figure}

\begin{figure}
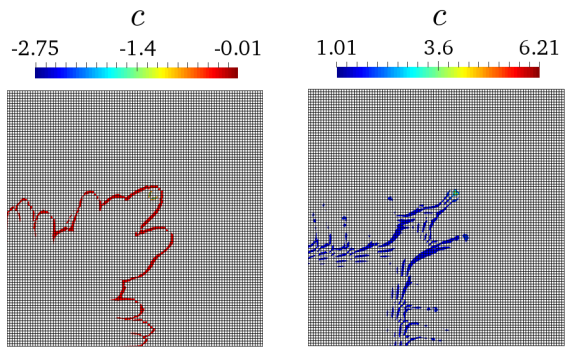

\centering 
  {\includegraphics[clip, scale=0.4]
  {Figures/conc_legend_low.png}}
  	\vspace{0.05 cm}
  	\hspace{0.43 cm}
  {\includegraphics[clip, scale=0.4]
  {Figures/conc_legend_up.png}} \\
  	\vspace{0.05 cm}
  %
  {\includegraphics[clip, scale=0.3]
  {Figures/supg_t_250_negative_val.png}}
  	\hspace{0.35 cm}
  %
  {\includegraphics[clip, scale=0.3]
  {Figures/supg_t_250_unphysical_val.png}}
  \caption{\textsf{Violation of maximum 
  principle:}~These figures show {\it unphysical} concentration values at time $t = 250$ 
  for the quarter five-spot problem under the SUPG stabilized formulation. 
  The parameters for these figures are the same as those used for Figure \ref{Fig:conc_case_I}.
  \label{Fig:conc_unphysical}}
\end{figure}

We performed numerical simulations using the SUPG formulation and its modifications, got by adding one or both of the isotropic and crosswind artificial diffusion terms to the SUPG formulation. We have used $100 \times 100$ quadrilateral finite elements. Table \ref{Tab:Input_case_II_dispersion} provides the parameters used in the numerical simulation. 
Figure \ref{Fig:conc_case_I} shows the concentration profiles at time levels $t = 100, 175$ and $250$ for three cases: the SUPG formulation, SUPG with isotropic diffusion, and SUPG with crosswind diffusion. The numerical solutions for the concentration field under all the three cases did not remain in the range $[0,1]$, and \emph{they violated the maximum principle and the non-negative constraint for all the three cases}. For instance, at $t = 250$, the concentration field under the SUPG formulation is $-2.75 \leq c \leq 6.21$; see Figure \ref{Fig:conc_unphysical}.

\begin{figure}
\centering 
  {\includegraphics[clip, scale=0.5]
  {Figures/conc_legend.png}} \\
  \subfigure[$t = 100$]
  {\includegraphics[clip, scale=0.3]
  {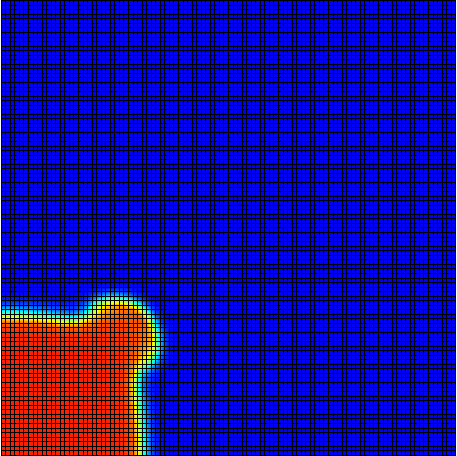}}
  \subfigure[$t = 175$]
  {\includegraphics[clip, scale=0.3]
  {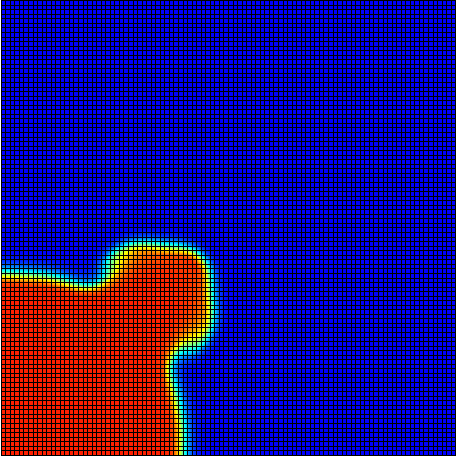}}
  \subfigure[$t = 250$]
  {\includegraphics[clip, scale=0.3]
  {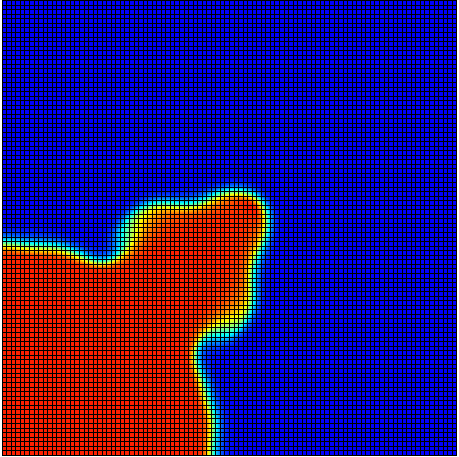}}
  \caption{\textsf{Suppression of physical instabilities:}~These figures show concentration profiles at various time levels for the quarter five-spot problem under the modified SUPG stabilized formulation with isotropic and crosswind artificial diffusion stabilizers. The parameters for these figures are the same as those used for Figure \ref{Fig:conc_case_I}. Although the results satisfy the maximum principle and the non-negative constraint for the concentration field, the utilized stabilizers suppress physical instabilities (i.e., the fingers).
  \label{Fig:conc_case_II}}
\end{figure}

Since adding a stabilizer one at a time---either isotropic diffusion or crosswind diffusion---did not eliminate spurious oscillations, we have solved the same problem by adding both the stabilizers to the SUPG formulation. Figure \ref{Fig:conc_case_II} shows the corresponding results for the concentration field. Although adding both the stabilizers have eliminated numerical instabilities (spurious oscillation, and the violations of the maximum principle and the non-negative constraint), it also suppressed physical hydrodynamic instabilities. \emph{That is, we did not observe VF under the SUPG formulation with isotropic and crosswind diffusion stabilizers.}

\section{CLOSURE}
\label{Sec:S5_Stab_CR}
We have shown that numerical stabilization techniques, which aim at removing spurious oscillations in the concentration and temperature profiles, may suppress physical instabilities (such as VF). Several popular stabilized formulations violate the non-negative constraint and the maximum principle for the concentration field. Several prior studies have reported such violations for the cases of advection-diffusion and VF with the single diffusivity under dispersion tensors with strong anisotropy (e.g., see \citep{chang2017variational}). However, we reported that such violations could occur even under isotropic diffusivity for VF with DD effects. 
 
The research reported in this paper highlights two points. First, practitioners who use existing numerical formulations should be wary of possible inaccurate and unphysical results when dealing with physical instabilities. They must choose with caution an appropriate numerical method that suppresses unphysical oscillations and captures physical instabilities at the same time.  Second, there is a need to develop new robust numerical formulations for coupled flow-transport problems, and in particular, to simulate double-diffusive effects. The non-negative formulations based on variational inequality, developed for coupled flow-transport problems with single diffusive effects \citep{chang2017variational} and flow through porous media problems \citep{mapakshi2018scalable}, could provide a viable route for developing numerical formulations for double-diffusive effects.

\bibliographystyle{plainnat}
\bibliography{Master_References/Master_References,Master_References/Books}  
\end{document}